\documentclass[twoside,12pt]{article}
\usepackage{amsmath,amssymb}
\date{}
\newtheorem{Lemma}{LEMMA}[section]
\newtheorem{Corollary}[Lemma]{COROLLARY}
\newtheorem{Theorem}[Lemma]{THEOREM}
\newtheorem{Proposition}[Lemma]{PROPOSITION}
\newtheorem{Definition}[Lemma]{DEFINITION}

\newcommand{\bnum}{\begin{enumerate}}
\newcommand{\enum}{\end{enumerate}}
\newcommand{\bi}{\begin{itemize}}
\newcommand{\ei}{\end{itemize}}
\newcommand{\btab}{\begin{tabular}}
\newcommand{\etab}{\end{tabular}}
\newcommand{\beq}{\begin{eqnarray*}}
\newcommand{\eeq}{\end{eqnarray*}}
\newcommand{\beqn}{\begin{eqnarray}}
\newcommand{\eeqn}{\end{eqnarray}}

\newcommand{\bq}{\begin{equation}}
\newcommand{\eq}{\end{equation}}

\newcommand{\CC}{{\cal C}}

\newcommand{\CL}{{\cal L}}
\newcommand{\CM}{{\cal M}}

\newcommand{\CR}{{\cal R}}

\newcommand{\CT}{{\cal T}}

\def\phi{\varphi}
\def\epsilon{\varepsilon}
\newcommand{\BS}{\mathbb S}
\newcommand{\BR}{\mathbb R}
\newcommand{\BC}{\mathbb C}

\newcommand{\BZ}{\mathbb Z}
\newcommand{\BH}{\mathbb H}
\newcommand{\DP}{{{\rm PG}(3,\BR)}}
\newcommand{\SO}{{{\rm SO}_3 \BR}}
\newcommand{\So}{{{\rm SO}_2 \BR}}
\newcommand{\oo}{{{\rm O}_2 \BR}}
\newcommand{\kasten}{\vbox{\hrule height 8pt width 8.6pt depth -7.4pt
    \hbox{\vrule width 0.6pt height 7.4pt
    \kern 7.4pt \vrule width 0.6pt height 7.4pt}
    \hrule height 0.6pt width 8.6pt}}
\newcommand{\ok}{\hfill\kasten}
\newcommand{\bpf}{\begin{Proof}}
\newcommand{\epf}{\ok\end{Proof}\bigskip\noindent}
\newcommand{\bthm}{\begin{Theorem}}
\newcommand{\ethm}{\end{Theorem}}
\newcommand{\ble}{\begin{Lemma}}
\newcommand{\ele}{\end{Lemma}}
\newcommand{\bprop}{\begin{Proposition}}
\newcommand{\eprop}{\end{Proposition}}
\newcommand{\bcor}{\begin{Corollary}}
\newcommand{\ecor}{\end{Corollary}}

\begin{document}
\title{Rotational spreads and rotational parallelisms and oriented parallelisms of ${\rm PG}(3,\BR)$}

\author{Rainer L\"owen}

\maketitle
\thispagestyle{empty}

\begin{abstract}
We introduce topological parallelisms of oriented lines (briefly called oriented parallelisms). Every topological 
parallelism (of lines) on $\DP$ gives rise to a parallelism of oriented lines, but we show that even the most homogeneous 
parallelisms of oriented lines other than the Clifford parallelism do not necessarily arise in this way. In fact we determine 
all parallelisms of both types that admit a reducible $\SO$-action (only the Clifford parallelism admits a larger group
\cite{dim4}), and it turns out surprisingly that there are far more oriented parallelisms of this kind than ordinary parallelisms.

More specifically, Betten and Riesinger \cite{nonreg} construct ordinary parallelisms by applying $\SO$ to rotational Betten spreads.
We show that these are the 
only ordinary parallelisms compatible with this group action, but also the `acentric' rotational spreads considered by them 
yield oriented parallelisms.
The automorphism group of the resulting (oriented or non-oriented) parallelisms is always $\SO$, no matter how large the 
automorphism group of the non-regular spread is. The isomorphism type of the parallelism 
depends not only on the isomorphism type of the spread used, but also on the rotation group applied to it.
We also study the rotational Betten spreads used in this construction and 
their automorphisms. 
 

\end{abstract}

\section{Introduction}

In \cite{nonreg}, Betten and Riesinger consider the so-called rotational Betten spreads on $\DP$, which were 
discovered by Betten \cite{bettenTE}. They admit
an action of $\So$ which induces axial collineations on the associated 4-dimensional translation plane. Betten and 
Riesinger observe that the 
Betten spreads  are composed of a family of reguli with common axis of symmetry and 
common center (a concentric family, as they call it), and that they are the only spreads that can be described in this way, 
see \cite{nonreg}, Section 4. 
They show that each of these spreads gives rise to a parallelism whose equivalence 
classes are the images of the given spread under a reducible $\SO$-action which extends the given $\So$-action. They also show that 
there are other families of reguli with common axis, called acentric, which also define spreads admitting an 
$\So$-action of the kind considered. 

We shall make the easy observation that the (concentric) Betten spreads are the only compact spreads in $\DP$ 
admitting an axial $\So$ action 
that extends to an action of $\oo = \So \cdot \BZ_2$, and we deduce that 
the parallelisms obtained by applying the reducible $\SO$-action to them are the only topological parallelisms that have this group
as their automorphism group, see Theorems \ref{charspread} and \ref{main-2} below. 
Thus, they are among the most homogeneous non-classical topological parallelisms. 
Indeed,  according to \cite{dim4} the 
`classical' Clifford parallelism is the only one 
with a group of dimension $\ge 4$, and the only other possible 3-dimensional group is $\SO$ with its action induced by the
irreducible $\rm Spin_3$ action on $\BR^4$. 
This follows from the fact that the automorphism group of a topological parallelism is compact \cite {zush}, \cite{unzush}.

We examine the proof given in \cite{nonreg} 
that the Betten spreads produce parallelisms with reducible $\SO$-action. It turns out that virtually 
the same proof can be applied to 
spreads defined by coaxial but acentric families of reguli. However, this time we do not obtain parallelisms of (non-oriented) 
lines as 
before, but parallelisms of oriented lines, see Theorem \ref{main-1}. 
Still their automorphism group is as large as possible for a non-classical parallelism. 
Again, according to Theorem \ref{main-2}, they are characterized by their group.

Therefore, we feel justified to define oriented parallelisms as a new kind of parallelisms. 
We record their basic properties, which resemble those 
of ordinary parallelisms except that their set of parallel classes is homeomorphic to the 2-sphere rather than the projective plane. 

As a by-product, our approach yields a very short proof of the result from \cite{cliff} that rotating a rotational 
spread produces the Clifford parallelism if and only if the spread is the regular (complex) spread and the rotation 
group fixes the center of the family of reguli, see Theorem \ref{cliff}.

Regarding spreads with an axial $\So$-action, we observe that they are precisely those obtained from coaxial (concentric or acentric) 
families of reguli. A precise description of the admissible parameters for such families of reguli has been obtained by 
Harald L\"owe (yet unpublished). We consider the automorphism group of such spreads and show that it is 1-dimensional except for
the regular spread and a small class of spreads with 2-dimensional groups. 
The latter spreads belong to some of the 4-dimensional translation planes with 7-dimensional group, which were classified by Betten
\cite{1fix}, \cite{2fix}, \cite{komstandgp}. To be more exact, only certain special instances of the planes described in 
Satz 1, 2 and 3 of \cite{2fix} occur here. In a majority of the cases, we describe these spreads in terms of families of reguli. 
This description is somewhat simpler than the original one. 

Finally we study isomorphisms and automorphisms of (oriented) parallelisms with reducible $\SO$-action. We show that the automorphism
group is always $\SO$,  even though the automorphism groups of the spreads differ in size, see Theorem \ref{autpar}. 
Correspondingly, the same spread rotated 
by different admissible copies of $\SO$ can yield non-isomorphic (oriented or non-oriented) parallelisms. 
In fact, typically this is what happens, see Theorem \ref{noniso}.

\section{Parallelisms versus oriented parallelisms}

An \it oriented line \rm of $\DP$ may be defined as a line together with a preferred direction of circuit or as a 
2-dimensional subspace of $\BR^4$ together with a vector space orientation. The spaces $\CL$ and $\CL^+$ of lines and of 
oriented lines are topologized as coset spaces of ${\rm PGL}(4,\BR)$. In this way, $\CL^+$ becomes a twofold covering space of
$\CL$. 

Suppose that $\CC \subseteq \CL$ is a \it compact spread \rm, i.e. a compact set of pairwise disjoint lines covering 
the point set of $\DP$. Then $\CC$ is homeomorphic to the 2-sphere $\BS_2$, see \cite{CPP}, 64.4(b), and hence its inverse 
image under the covering map $\CL^+ \to \CL$ is a disjoint union of two 2-spheres, each of which is mapped onto 
$\CC$ and hence is a compact spread consisting of 
oriented lines. We call each of these spheres an \it oriented spread. \rm Thus we have proved

\bprop\label{orspread}
Every compact spread of $\DP$ can be oriented in exactly two ways.
\hfill \kasten
\eprop

A \it topological parallelism \rm on $\DP$ is a set $\Pi$ of pairwise disjoint compact spreads covering $\CL$ which is a compact
subset of the hyperspace of all compact subsets of $\CL$ (endowed with the Hausdorff topology). Similarly, an \it oriented 
topological parallelism \rm is a compact set $\Pi^+$ of pairwise disjoint compact, oriented spreads covering $\CL^+$. 
From \ref{orspread} we deduce

\bprop\label{orpar}
Every topological parallelism $\Pi$ of $\DP$ gives rise to a unique oriented topological parallelism $\Pi^+$, which consists of 
all oriented spreads covering members of $\Pi$.
\hfill \kasten
\eprop

As we shall see, there is no converse to this proposition. \\

For the sake of distinction, we shall sometimes refer to non-oriented parallelisms as \it ordinary parallelisms. \rm The properties
of the two kinds of parallelisms are very similar. See, e.g., the section `Preliminaries' of \cite{unzush} for a brief account 
of the essentails in the ordinary case. One distinctive feature is the homeomorphism type of a parallelism, considered as a 
subset of the hyperspace of $\CL$ or $\CL^+$:

\bprop\label{topofpi}
An ordinary topological parallelism $\Pi$ is homeomorphic to the real projective plane. An oriented topological parallelism 
$\Pi^+$ is homeomorphic to the 2-sphere. 
\eprop

\bpf Choose a point $p$ and let $\CL_p \subseteq \CL$ and $\CL_p^+\subseteq \CL^+$ be the stars of non-oriented lines 
and of oriented lines passing through $p$, respectively. The ordinary star is a projective plane, and the oriented one is a 
2-sphere. The map sending a parallel class to its unique member in the 
star is a homeomorphism in both cases.
\epf

The last result enables a nice way of proving the well known fact that the space $\CL^+$ of oriented lines of $\DP$ is 
homeomorphic to the product $\BS_2\times \BS_2$ of two 2-spheres. One classical proof of this property uses Pl\"ucker 
coordinates, which provide a homeomorphism $\CL \to Q$ onto the Klein quadric of index 3 in ${\rm PG}(5,\BR)$. It is easy to see that 
$Q$ is homeomorphic to the quotient of $\BS_2 \times \BS_2$ with the identification $(x,y) = (-x,-y)$; 
a detailed proof is given in \cite{dim3}. From this it follows
that the twofold covering space $\CL^+$ is $\BS_2 \times \BS_2$. 

Another proof is provided in \cite{kling}, p. 290 -- 291, and attributed to E. Study. 
Identify $\BR^4$ with the quaternion skew field $\BH$. Let $\Lambda$ and $\Theta$ be the subgroups of 
${\rm SO}(4,\BR)$ consisting of the maps $q\mapsto aq$ and $q \mapsto qb$, respectively, where $a,b$ are quaternions of norm 1.
Then an oriented line $L$ is mapped to the pair of its stabilizers $(\Lambda_L,\Theta_L)$, both endowed  with orientations 
matching the orientation of $L$. Now $\Lambda_L$ and $\Theta_L$ are both isomorphic to $\So$. Together with the orientation, 
each of them defines a unique pure quaternion of norm 1, that is, an element of $\BS_2$. For instance, $\Lambda_L$ contains 
exactly two elements of order 4. One of them, say $\lambda$, comes first when $\Lambda_L$ is traversed in the sense 
of orientation, starting from  the identity. The pure quaternion $a$ associated to $\Lambda_L$ is the one satisfying 
$\lambda(q) = aq$.

The following proof also uses the groups $\Lambda$ and $\Theta$; the orbits of $\Lambda$ and $\Theta$ on $\CL^+$ form the
left and right oriented Clifford parallelism on $\DP$, respectively.

\bprop
Let $\Pi^+_\Lambda$ and $\Pi^+_\Theta$ be the left and right oriented Clifford parallelisms of $\DP$, as above. Then the map
   $$\alpha: \CL^+\to  \Pi^+_\Lambda\times \Pi^+_\Theta \approx \BS_2 \times \BS_2 $$
sending an oriented line $L$ to the pair $(\Pi^+_\Lambda(L), \Pi^+_\Theta(L))$ of its parallel classes is a homeomorphism.   
\eprop

\bpf
The map is surjective because it is equivariant with respect to the transitive ${\rm SO}(4,\BR)$-actions on both spaces.
Injectivity of $\alpha$ follows from transitivity together with the observation that the stabilizer
of $L$ coincides with that of $\alpha (L)$. 

Alternatively, we show that $\alpha$ is injective by proving that
   $$\Pi_\Lambda^+(L) \cap \Pi_\Theta^+(L) = \{L\}.$$  
By transitivity, we may assume that $L$ is the 2-dimensional subspace $\langle 1, q \rangle^+ \le
\BH$ with oriented basis $\{1,q\}$, where $q$ is a pure quaternion of norm 1. Suppose that there are
quaternions  $a$, $b$ of norm 1 such that the images of $L$ under $x \mapsto ax$ and $x \mapsto xb$ are equal 
(as oriented lines), that is,
   $$ \langle a,aq\rangle^+ = \langle b, qb\rangle^+ =M. $$
We have to show that $M =L$ as oriented lines.   
The intersection of $M$ with the sphere of radius 1 is an orbit of a one-parameter group of $\Theta$, hence there is a quaternion
$c$ of norm 1 such that the map $\vartheta = (x\to xc)\in \Theta$ fixes $M$ and sends the vector $b \in M$ to $a \in M$. 
Then the above equation
implies that 
             $$\langle a,aq\rangle^+ = M = \vartheta(M) = \langle a, qa\rangle^+.$$ 
Since there is only one oriented orthonormal basis for $M$ containing  $a$ as its first element, we conclude that $aq = qa$, 
i.e., that $a$ belongs to 
the centralizer of $q$ in $\BH^\times$, which equals $L\setminus\{0\}$ and is a subgroup of $\BH^\times$. 
It follows that $M=L$.           
\epf

\bf Remark. \rm In the sequel, we shall never consider non-topological spreads or parallelisms. Therefore, we shall usually omit 
the words `topological' or `compact' (in the case of spreads) unless we want to stress that we have just constructed a 
topological object.

\section{Automorphisms}\label{aut}

The automorphism group $\Sigma$ of an ordinary or oriented parallelism $\Pi$ or $\Pi^+$ is defined as the group of 
all collineations of $\DP$
that preserve $\Pi$ or $\Pi^+$, respectively. In the ordinary case, we proved that this group is compact \cite{unzush}, \cite{zush}. 
The same proof goes through in the oriented case virtually unchanged. Lemma 3.4 of \cite{unzush} states that every automorphism 
$\sigma$ induces 
equivalent actions on the line stars of any two fixed points. Of course, this has to be adapted by using the oriented line stars.
The applications of this lemma work just as well. In the proof of compactness, the goal is to rule out all types of unbounded 
cyclic subgroups of $\Sigma$. The contradiction always arises in the form that certain lines are shown to be be parallel by a 
continuity argument, but they are distinct and have a common point. This contradiction does not depend on the orientations, 
hence it is not even necessary to keep track of orientations in order to adapt the proof. Thus we have

\bthm\label{comp}
Like in the ordinary case, the automorphism group of an oriented topological parallelism $\Pi^+$ on $\DP$ is compact.
\hfill \kasten
\ethm

The automorphism groups of an ordinary parallelism and of its associated oriented parallelism are the same, so by 
\cite{coll} we get part (a) of the next theorem; compare also \cite{unzush}, Corollary 1.2.
The proof of part (b) is again virtually the same as in the oriented case, see \cite{dim4}. 

\bthm a) The automorphism group of oriented Clifford parallelism is the 6-dimen\-sional group ${\rm PSO}(4,\BR) \cong 
\SO \times \SO$. 

b) A non-classical oriented topological parallelism has an automorphism group of dimension at most 3.
\hfill \kasten
\ethm

\section{Rotational spreads of $\DP$}\label{rot}

For general information on compact spreads and their relationship to topological 
translation planes, we refer the reader to \cite{CPP}, Section 64.
The following considerations are inspired by \cite{nonreg}. 
We assume in this section that we have a (compact) spread $\CC$ of $\DP$ that admits an axial action of $\So$. When we consider 
$\CC$ as a fibration of $\BR^4\setminus \{0\}$ by 2-dimensional subspaces, this means that $\Phi \cong \So$ acts trivially on some  
2-space $S \in \CC$ and induces the ordinary $\So$ on some complementary subspace $W$. 

We prefer, however, to represent the spread as a
set of lines of $\DP$; 
then $\Phi$ fixes one line $Z$ pointwise and fixes one other line $V$ disjoint from $Z$. The line $Z$ belongs to $\CC$ by assumption,
but in fact this can be proved, see \ref {charspread} below.
Moreover, $\Phi$ fixes all hyperplanes spanned by $V$ together with a point of $Z$. We choose one of these hyperplanes, $F$,
and consider the affine complement $\BR^3$ of $F$, endowed with a $\Phi$-invariant inner product. We may assume that $Z$ is 
the $z$-axis $Z$, spanned by the third standard basis vector $e_3$, and $V$ is the line at infinity of the plane 
$\langle e_1, e_2\rangle$. Thus, the group $\Phi$ consists of the ordinary rotations with axis $Z$. Its orbits on the line space 
are two fixed lines, $Z$ and $V$, and all reguli carried by
$\Phi$-invariant one-sheeted hyperboloids. Since $\CC$ is a 2-sphere, there must be a second fixed line in $\CC$, thus
$V \in \CC$. The reguli belonging to $\CC$ have to be either all right-screwed or all left-screwed, or else $\CC$ 
could not be connected. Right-screwed reguli and left-scewed ones are exchanged by the map $(x,y,z) \to (x,-y,z)$, hence
we shall not pay attention to the choice between the two.

We see that $\CC$ can be described by specifying the one-sheeted hyperboloids carrying the 
reguli contained in $\CC$. This can be done by identifying the hyperbola branches obtained by intersecting the 
hyperboloids with the half plane $E = \{ (x,0,z)\ \vert \ x,z \in \BR, x>0\}$. Such a hyperbola is determined by three
parameters, namely the coordinates of its vertex $(r,0,b)$ (the point closest to $Z$), and the slope of its upper asymptote, $a>0$.
By the defining property of a spread, any two of the hyperbolae are disjoint. Hence they have distinct $r$-parameters, and we may 
label the hyperbolae as $H_r$, $0 < r \in \BR$. In fact, every positive $r$ occurs, because the $H_r$ must cover $E$. 
Thus, $b$ and $a$ become functions of $r$. The function $a$ is  decreasing, because the hyperbolae are pairwise disjoint, and in 
fact strictly decreasing, because two hyperboloids with the same $a$ have the same set of points at infinity. For $r \to 0$ and 
$r \to \infty$, the reguli must converge to $Z$ and to $V$, respectively, which means that $a(r) \to \infty$ in the first case and 
$a(r) \to 0$ in the second. 

The elements of a spread $\CC$ as considered here (except $Z$ and $V$) are determined 
by two parameters $r>0$ and $\phi\in \Phi$, 
hence we see that $\CC$ is compact (in fact, a 2-sphere).

We now distinguish two cases; 

\it Case 1 (the concentric case): \rm the function $b$ is constant. 

\it Case 2 (the acentric case): \rm the function $b$ is not constant. 

In Case 1, the conditions on the function $a$ stated above
are sufficient to ensure that the hyperbolae $H_r$ are pairwise disjoint and cover $E$. The points $p \in F$ at infinity are also 
simply covered by $\CC$, and we obtain a spread, which is called a \it rotational Betten spread \rm 
in \cite {nonreg}. These spreads include 
the classical (regular, complex) spread, which is obtained precisely when $a(r) = {a(1)\over r}$, see Theorem \ref{regular}.
In Case 2, extra conditions on the function $b$ are needed to ensure the disjointness of the hyperbola branches. 
The conditions have been determined by Harald L\"owe (yet unpublished). Examples of spreads of this type are given in \cite{nonreg}.
Definitely, the class of spreads in Case 2 is much larger than the class defined by Case 1.

It is obvious that Case 1 occurs if and only if the 
spread $\CC$ is invariant under the orthogonal reflection $\sigma$ in some line perpendicular to $Z$, so that $\CC$ admits the 
group $\Psi = \Phi \cdot \langle \sigma \rangle \cong \oo$. 

\begin{Definition} \rm We shall call the spreads described here (as well as all isomorphic copies) \it $\So$-rotational spreads\rm,
and the special ones arising in Case 1 will be called \it $\oo$-rotational spreads. \rm The names are justified by the following
\end{Definition}

\bthm\label{charspread}
Let $\Gamma$ be either  $\oo$ or $\So$. 
A (compact) spread $\CC$ of $\DP$ is a $\Gamma$-rotational spread if and only if it is invariant under 
an action of the group $\Gamma$ on $\DP$ that is equivalent to the projective extension of the standard action of 
$\Gamma \le \SO$ on $\BR^3$.
\ethm

\bpf The `if' part being obvious, we merely consider  the `only if' assertion.
It suffices to show that $\CC$ satisfies the basic assumption of the preceding considerations, namely, that the action of the 
identity component $\Phi \le \Gamma$ is axial. By the hypothesis of the theorem, $\Phi$ fixes some line $Z$ pointwise;
we have to show that $Z \in \CC$.
For the corresponding fibration of $\BR^4$ and  the lifted action of $\Phi$, we see that some 2-space $S$ is fixed elementwise. 
It follows that $S$ belongs to the fibration, or 
else $S$ would be a Baer subplane of the translation plane defined by the fibration, 
and there would exist at most  one non-trivial automorphism fixing $S$ pointwise, see \cite{CPP}, 55.21.
\epf

\bthm\label{regular}
Let $\CC$ be a $\Gamma$-rotational spread defined by the functions $a(r)$ and $b(r)$. 
The spread $\CC$ is regular if and only if $b(r)$ is constant and $a(r) = {a(1) \over r}$.
In this case, the group $\mathop{\rm Aut}\CC$ of all collineations preserving $\CC$ 
is known to be 7-dimensional. 
\ethm

\bpf
If $b$ is not constant, then $\CC$ is not an $\oo$-spread and certainly not classical. If $b$ is constant, we may assume 
that $b \equiv 0$, and then the
assertion is proved in \cite{bettenTE}; note that Betten's parameter $t$ corresponds to our $1 \over r$, and his function $F(t)$
is our $a({1\over r})$; compare \cite{nonreg}. The automorphism group of the regular spread is $(\mathop{\rm GL_2}\BC)/\BR^\times$.
\epf

Now we consider the non-regular cases.

\bprop\label{charsubgp}
Let $\CC$ be a non-regular $\Gamma$-rotational spread. Then the automorphism group 
$\Sigma = \mathop{\rm Aut}\CC \le {\rm PGL}_4\BR$ is at most 2-dimensional, and 
the identity component $\Phi = \Gamma^1 $ is the only subgroup isomorphic 
to $\So$ in $\Sigma$. In particular, $\Phi$ is a characteristic subgroup of $\Sigma$.
\eprop

\bpf
If $\Sigma$ contains a 2-torus $\Delta$, then $\CC$ is regular. Indeed, all 2-tori in ${\rm PGL}_4\BR$ are conjugate, 
hence we may assume that the 2-torus 
is induced by the group of complex diagonal matrices with entries of norm one. This group is ineffective on the 
fibration of $\BR^4$ corresponding to $\CC$. The one-dimensional kernel of ineffectivity contains no elements with 
eigenvalue 1 except the identity because the line at infinity of the associated translation plane is fixed, 
hence it acts like the group of complex scalar 
matrices, and the claim follows.

If  $\Sigma$
contains several copies of $\So$ but not a 2-torus, then $\dim \Sigma \ge 3$, and the translation plane defined by $\CC$ has an 
automorphism group of dimension at least $4 + 3  + 1 = 8$. In fact, the latter dimension equals 8 since we assumed that 
$\CC$ is non-regular, see \cite{CPP}, 73.1. According to \cite{CPP}, 73.10, 73.11, there is a unique 
4-dimensional translation plane with a solvable 
8-dimensional automorphism group $\Xi$. The stabilizer $\Xi_0$ of the origin acts on $\BR^4$ 
by lower triangular matrices, and  does not contain a non-trivial torus group. 
Similarly, the translation planes with a non-solvable automorphism group do not admit an axial torus group, see \cite{CPP}, 73.13 and 
73.19 combined with 73.18(e). Thus we end up with a contradiction and the theorem is proved.
\epf

\bcor\label{isom}
If $\CC$ and $\CC'$ are non-regular $\Gamma$-rotational spreads and $\Phi \cong \So$ is contained in the 
automorphism groups of both spreads, then every isomorphism $\xi:  \CC \to \CC'$ belongs to the normalizer $\mathop{\rm Ns}\Phi$
in ${\rm PGL}_4\BR$. \hfill \kasten

\ecor 

The next theorem shows that rotational spreads with  $\dim \Sigma = 2$ are determined by finitely many real parameters. 
Thus, a generic $\Gamma$-rotational spread has a 
one-dimensional automorphism group.

\bthm\label{dim2}
Let $\CC$ be a non-regular $\Gamma$-rotational spread with a 2-dimensional automorphism group. Then the translation plane 
defined by $\CC$ is one of the following planes described in \cite{2fix}:

(1) If $\CC$ is $\oo$-rotational, then the plane is one of the planes of Satz 1, and the parameters $p,c,w$ 
satisfy $p = c = 0$, with no further restriction on $w \in \left ]0,1\right [$.

(2) If $\CC$ is not $\oo$-rotational, then the plane is either

(2a) any of the remaining planes of Satz 1 with parameters $p= 0\ne c$, or

(2b) any of the planes of Satz 2 with $p=0$ and $\vert d \vert \ge {1\over 2}$, or

(2c) any of the planes of Satz 3 with parameters $p,q,c,d$ satisfying $p = 0 < q$ and $d > 0$, as well as the conditions (2) stated in \cite{2fix}.
\ethm

The properties of the planes are described exhaustively in \cite{2fix}. See also \cite {knarr}, p. 77, 
where some minor details about the isomorphism types of the planes given in \cite{2fix}, Satz 3 are corrected.

\bpf 
The 4-dimensional translation planes with a 7-dimensional automorphism group have been classified by Betten, 
\cite{1fix}, \cite{2fix}, \cite {komstandgp}, and complete information on the spreads and the automorphisms is given there. 
The planes of \cite{1fix} and of \cite {komstandgp}  do not admit any torus action.
Inspection of the planes of \cite{2fix} shows that only the planes identified in the theorem admit 
an axial $\So$-action, and only those in (1) admit $\oo$.
\epf

We close this section by describing some of the $\Gamma$-rotational spreads with 2-dimensional automorphism group in terms of
the functions $a(r)$ and $b(r)$ that define the hyperboloids whose right-screwed reguli make up the spread. We do this in the cases 
(1), (2a) and (2b) of Theorem \ref{dim2}. In principle, it can also be done in case (2c), but the result is not too pleasant.

Following the notation of Betten \cite{2fix}, we let $\BR^4 = W \oplus S$, and all spreads are given as 
fibrations of $\BR^4 \setminus \{0\}$ by 2-dimensional subspaces. All fibrations shall  contain $W$ and $S$ as elements, 
and $S$ is the axis of rotation.  
The remaining 2-dimensional subspaces are given in \cite{2fix} in the form $\{(v,Av) \ \vert \ v\in W\}$, where 
$A$ is a regular $2\times 2$ matrix.
We shall select some of these matrices so that we get a cross section to the $\So$-orbits in the fibration. 
This is achieved by setting the rotation parameter $\varphi$ equal to zero. Then we pass to 
affine coordinates in $\BR^3$ by setting the third coordinate equal to 1, and compute the distance $r$ of the resulting line from the 
$z$-axis as well as the $z$-coordinate $b(r)$ of the closest point and the slope $a(r)$ of the given line.

For the planes of \cite{2fix}, Satz 1 with $p=0$, the cross section is given by the matrices
  $$A(s) =\left ( \begin{matrix}  s & 0 \cr
                               s^wc & s^w \cr 
                  \end{matrix} \right ), \quad s>0.$$
Here, $w \in \left ] 0,1\right [$ and $c \in \BR$ are parameters defining the fibration, 
and $s >0$ parametizes the cross section.                   
Taking points of the 2-dimensional subspace defined by $A(s)$ and dividing by the third coordinate, we obtain the points 
of the corresponding line in $\BR^3$ as
   $$({1\over s}, {y\over sx}, s^{w-1}(c + {y \over x})).$$
The point closest to the $z$-axis is obtained by setting $y = 0$, hence $r = s^{-1}$, and the spreads arising from \cite{2fix} 
Satz 1 with $p=0$ are described by        
    $$a(r) = r^{-w}, \ \ b(r) = r^{1-w}c.$$
In the special case $c = 0$ we obtain the $\oo$-rotational spreads with 2-dimensional automorphism group.         

The planes of \cite{2fix}, Satz 2 with parameter $p = 0$ are given by
            $$A(t) =\left ( \begin{matrix}  e^t & 0 \cr
                                         te^t & de^t \cr 
                  \end{matrix} \right ), \quad t \in \BR.$$
The parameter $d$ satisfying $\vert d \vert \ge {1\over 2}$ defines the fibration. Proceeding as before, 
for fixed $t$ we obtain a line in $\BR^3$ with points
            $$(e^{-t}, e^{-t}{y \over x}, t + d{y\over x}).$$
We see that $r = e^{-t}$, and the planes of \cite{2fix}, Satz 2 with $p = 0$ are described by
        $$a(r) = {d\over r},\ \  b(r) = -\ln r.$$            
Recall from Theorem \ref{regular} that $a(r) = {d\over r}$, $b(r) = 0$ defines the regular complex spread.

\section{$\SO$-invariant parallelisms of $\DP$ obtained from rotational spreads}

The group $\SO$ can act on $\DP$ either irreducibly or reducibly. The lifted action on $\BR^4$ is equivalent to the action
of the group ${\rm Spin}_3$ of quaternions of norm 1 on the quaternion skew field $\BH$ given by $q \to aq$ or 
$q \to aqa^{-1}$, respectively, where 
$q \in \BH$ and $a \in {\rm Spin}_3$. The latter, reducible action has an affine description on $\BR^3$ as 
the ordinary rotation group $\Omega$. 

\begin{Definition} \rm Let $\mathop{\rm AGL}_3\BR$ be the group of affine transformations of $\BR^3$, and
let $\Omega\le \mathop{\rm AGL}_3\BR$ be any affine  group isomorphic to $\SO$.  Let 
$\Gamma \le \Omega$ be a subgroup isomorphic to $\So$ or to $\oo$, and let $\CC$ be a spread invariant under the natural action of 
$\Gamma$ (so $\CC$ is a $\Gamma$-rotational spread by Theorem \ref{charspread}). 
Then we say that $\Omega$ is \it $\Gamma$-admissible \rm for $\CC$. Note that for $\Gamma = \oo$ this implies that 
$b(r)$ is constant and that $\SO$ fixes the common center of the reguli. Up to conjugacy in $\mathop{\rm AGL}_3\BR$,
we may always assume that $\Phi = \Gamma^1$ consists of the ordinary rotations about the $z$-axis $Z$.
\end{Definition}

\ble\label{admiss}
Fix a group $\Gamma \le \mathop{\rm GL}_3\BR$ isomorpic to $\So$ or to $\oo$ 
such that $\Phi = \Gamma^1$ consists of the ordinary rotations about the $z$-axis $Z$.

If a rotation group $\Omega\le \mathop{\rm AGL}_3\BR$ is $\Gamma$-admissible, then all 
other $\Gamma$-admissible rotation groups $\Omega^\prime\le \mathop{\rm AGL}_3\BR$ are obtained as conjugates 
$\Omega^\xi$, where $\xi(x,y,z) = (x,y,sz+t)$ with $s,t \in \BR$ and $s>0$; if $\Gamma = \oo$, then only $t=0$ is allowed. 
\ele

\bpf
All rotation groups are conjugate in $\mathop{\rm AGL}_3\BR$, and all one-parameter groups of $\Omega^\prime$ 
are conjugate in $\Omega^\prime$.
Therefore, $\Omega^\prime = \Omega^\alpha$ for some affine map $\alpha$ normalizing $\Phi = \Gamma^1$. Thus $\alpha$ 
fixes the rotation axis $Z$ of $\Phi$ and the orthogonal space $Z^\perp$ that is invariant under $\Phi$. 
The map induced by $\alpha$ on $Z$ must be of the form $z \mapsto sz+t$, and on $Z^\perp$ a homothety is induced. 
Now the homotheties of $\BR^3$ centralize $\Omega$, and the assertion for $\Gamma = \Phi$ follows. If $\Gamma = \oo$, then 
both $\Omega$ and $\Omega^\prime$ fix the unique fixed point of $\Gamma$, whence $t = 0$.
\epf

\begin{Definition} \rm If $\CC$ is $\Gamma$-rotational and $\Omega$ is $\Gamma$-admissible, we define a set of spreads by
   $$\Omega (\CC) = \{\omega(\CC) \ \vert \ \omega \in \Omega\}.$$
Furthermore, let $\CC^+$ be one of the two oriented spreads obtained from $\CC$ and define $\Omega (\CC^+)$ in the same manner. 
\end{Definition}
 
Our first main result is the following theorem. Assertion (1) is contained in \cite {nonreg} under the hypothesis 
that $\CC$ is a rotational Betten 
spread as defined earlier.

\bthm\label{main-1}
Let $\CC$ be a $\Gamma$-rotational spread and let $\Omega$ be a $\Gamma$-admissible affine rotation group. \\
(1) If $\Gamma = \oo$, then the set $\Omega (\CC)$ of spreads is a topological parallelism of $\DP$. \\
(2) If $\CC$ is  
not an $\oo$-spread or $\Omega$ is not $\oo$-admissible, then $\Omega (\CC)$ is not a parallelism, 
but $\Omega (\CC^+)$ is an oriented 
topological parallelism. \\
(3) The oriented parallelism in assertion (2) does not arise from any ordinary parallelism.
\ethm

More explicitly, the assumption of part (2) means that either the function $b(r)$ is not constant or $b(r) \equiv b$ 
but $\Omega$ does not fix the point $(0,0,b))$.

\bpf
Assertion (3) is an immediate consequence of (2), because the image of $\Omega(\CC^+)$ under the coverng map 
$\CL^+\to \CL$ is $\Omega(\CC)$. 
The negative part of assertion (2) follows from 
the fact that  a topological parallelism $\Pi$ is homeomorphic to the real projective plane; therefore, 
the stabilizer of a transitive 
$\Omega$-action on $\Pi$ is $\oo$. Hence, the spreads $\CC\in \Pi$ are $\oo$-spreads, 
and $\Omega$ is $\oo$-admissible. We shall now prove (1) 
and the positive part of (2) simultaneously. We shall use $\CL^*$ as a shorthand for $\CL$ or $\CL^+$, whichever 
is appropriate, and similarly for other sets of lines. Note that, in general, $\CM^+$ denotes the full inverse image of
$\CM \subseteq \CL$ under the covering map $\CL^+ \to \CL$, whereas $\CC^+$ denotes just one of the two connected 
components of the inverse image. 

We may assume that $\Omega$ fixes the origin $0\in \BR ^3$ and that $\CC$ is constructed as in Section \ref{rot}, 
with respect to the group $\Phi$ of
rotations about the $z$-axis $Z$; compare Theorem \ref{charspread}.
First we describe the line orbits of $\Omega$. Clearly, $\Omega$ is 
transitive on the sets $\CL_0^*$ and $\CL_F^*$,
where the subscript $0$ refers to lines passing through the origin and the subscript $F$ indicates lines contained in the 
hyperplane $F$ at infinity. The remaining orbits are the sets $\CT_d^*$ consisting of the lines or oriented lines 
tangent to the sphere
of radius $d>0$ centered at $0$. We have to show two things, namely\\

(i) $\CC^*$ contains lines or oriented lines from every $\Omega$-orbit, and\\

(ii) If a line $L \in \CC^*$ has an image $\omega(L) \in \CC^*$ where $\omega \in \Omega$, then  $\omega \in \Gamma$. \\

Indeed, (i) expresses that $\Omega (\CC^*)$ covers $\CL^*$, and
(ii) expresses that any two $\omega$-images of $\CC^*$ are either disjoint or equal. 
We shall use the 
notation of Section \ref{rot}.
Condition (i) is obtained as follows. For the two special orbits  the assertion is obvious, since $Z\in \CL_0^*$ and $V\in \CL_F^*$, 
no matter which orientation we choose for $Z$ and $V$. Now the half plane $E$ is simply covered by the hyperbola 
branches $H_r$, and every 
$H_r$ contains a unique point closest to the origin, say at distance $d = d(r)$. Then the regulus $\CR_r \subseteq \CC$ 
carried by the hyperboloid 
corresponding to 
$H_r$ is contained in $\CT_d$. Since the hyperbolae are pairwise disjoint and cover $E$, every value $d$ occurs exactly once.
This proves (i). 

Condition (ii) is easily checked for the two $\Gamma$-invariant lines $Z$ and $V$ (with or without orientation). 
Now suppose that $L\in \CR_r^* \subseteq \CC^*$ and $\omega(L) \in \CC^*$. Then $\omega(L) \in \CT_{d(r)}^*$ because rotations leave
$\CT_d^*$ invariant, and $\omega(L)\in \CR_r^*$ by injectivity of the map $r \to d(r)$. Therefore, there is a rotation 
$\phi \in \Phi \le \Gamma$ such that $\phi\omega(L) = L$ as nonoriented lines. If $\Gamma = \oo$, then it follows that 
$\phi\omega \in \Gamma$ and, hence, $\omega \in \Gamma$. If $\Gamma = \So$, then $\phi\omega(L)\in \CC^+$ 
implies that $\phi\omega(L) = L$ as oriented lines, and then $\phi\omega = \rm id$. This proves (ii). 
Finally, compactness of the group $\Omega$ implies that $\Omega(\CC^*)$ is compact and, hence, is a topological 
parallelism or oriented parallelism. 
\epf

\bthm\label{main-2}
Every topological parallelism or oriented parallelism of $\DP$ admitting a reducible $\SO$-action is isomorphic to one of the 
examples given by Theorem \ref{main-1}.
\ethm

\bpf 
Up to conjugacy in $\mathop{\rm PGL_4\BR}$, there is only one reducible action of $\SO$ on $\DP$. Thus we may assume that $\Omega$
is $\SO$ acting on the affine space $\BR^3$ in the ordinary way. 
Now let $\Pi$ be a parallelism invariant under this action, possibly an oriented one. Then $\Omega$ does not 
act trivially on $\Pi$, because the line orbits of 
$\Omega$ are not contained in spreads. Thus $\Omega$ is transitive on $\Pi$, because $\Omega$ does not contain any 2-dimensional
closed subgroups. We know that $\Pi$ is a 2-sphere or a projective plane according as $\Pi$ is oriented or not, 
see  Proposition \ref{topofpi}, and hence 
the stabilizer $\Gamma$ of a spread
$\CC^* \in \Pi$ is $\So$ or $\oo$, respectively. According to Theorem \ref{charspread}, $\CC$ is a $\Gamma$-rotational spread.
Moreover, $\Omega$ is $\Gamma$-admissible for $\CC$, and the theorem follows.\epf

It is tempting to call the parallelisms characterized by Theorem \ref{main-2} \it rotational parallelisms. \rm However, 
this name has been given to a different class of parallelisms (the 3-dimensional regular parallelisms with 
2-dimensional automorphism group), see \cite{dim3}.

Let us now consider the regular case. The following result has been proved in \cite{cliff}, Theorem 8.1b in two distinct, 
very explicit ways, but the possible effect of changing the group $\Omega$ is not considered. Here we give a short conceptual proof. 

\bthm\label{cliff} Suppose that $\CC$ is the complex, regular spread and that $\Omega = \SO$ is $\oo$-admissible for $\CC$.
Then  $\Omega(\CC)$ is the Clifford parallelism. \ethm

\bcor
Clifford parallelism is the only regular topological parallelism 
admitting a reducible $\SO$-action.
\ecor

\bpf
Assume the hypothesis of \ref{cliff}. According to Theorem \ref{main-2}, Clifford parallelism is of the form $\Omega'(\CC')$. 
There are axial subgroups $\Phi \le \Omega$ and 
$\Phi^\prime \le \Omega^\prime$ isomorphic to $\So$. Let $Z^{(\prime)}\in \CC^{(\prime)}$ be the axis of $\Phi^{(\prime)}$ 
and $V^{(\prime)}\in \CC^{(\prime)}$ the 
second fixed line of $\Phi^{(\prime)}$. The group $\Sigma = \mathop{\rm Aut}\CC$ is doubly transitive on $\CC$, hence there is a collineation 
$\alpha$ of $\DP$ that sends $\CC^\prime$ to $\CC$ and $(Z^\prime, V^\prime)$ to $(Z,V)$. Replacing $\Omega^\prime$ with 
the conjugate $\Omega^{\prime \alpha}$ we obtain that $\CC = \CC^\prime$, $Z=Z^\prime$ and $V = V^\prime$. Now $\Omega^{(\prime)}$ fixes some hyperplane
$F^{(\prime)}$ containing $V$, and since $\mathop{\rm Ns}_\Sigma\Phi$ is transitive on the set of such hyperplanes we may assume that $F = F^\prime$.
Finally, the subgroups $\oo \le \Omega^{(\prime)}$ containing $\Phi$ can be adjusted to be equal using an affine translation in the direction of $Z$.

We have now achieved 
that the groups $\Omega$ and $\Omega'$ are both affine with respect to the same affine space and are both $\oo$-admissible for $\CC$. 
By Lemma \ref{admiss}, this implies that $\Omega' = \Omega^\xi = \xi^{-1}\Omega \xi$,
where the linear map $\xi$ is given in suitable coordinates by $(x,y,z)\mapsto (x,y,sz)$, for some $s >0$. Now $\rho: (x,y,z) \mapsto (sx,sy,z)$ is an 
automorphism of $\CC$, and $\xi\rho: (x,y,z) \mapsto (sx,sy,sz)$ centralizes $\Omega$. Thus we have
   $$\Omega(\CC) = \Omega^{\xi\rho}(\CC) = \rho^{-1}\Omega^\xi(\rho\CC) = \rho^{-1}\Omega'(\CC) \cong \Omega'(\CC). $$
This proves  Theorem \ref{cliff}, and the corollary follows in view of Theorem \ref{main-2}.   
\epf

Note, however, that the classical spread can yield different regular oriented parallelisms, 
if a rotation group is applied which does not fix the common center of the reguli (Theorem \ref{main-1} (2)). Next, we look at the non-regular case,
where in contrast to the previous result it will turn out that the choice of an admissible group $\Omega$ does matter, see Theorem \ref{noniso}.

\bthm\label{autpar}
Let $\CC^*$ be a (possibly regular) $\Gamma$-rotational spread (with orientation if $\Gamma = \So$), 
and let $\Omega = \SO$ be $\Gamma$-admissible for $\CC$. Suppose that the (oriented or
non-oriented) parallelism $\Pi = \Omega(\CC^*)$ is not Clifford. 
Then its full automorphism group $\Sigma = \mathop{\rm Aut} \Pi$  is equal to the group $\Omega = \SO$ 
that was used to construct it. 
\ethm

Note that the theorem applies to the oriented parallelisms obtained from the complex spread by rotating it about a point other than 
the center of the family of reguli.

\bpf 
We know that $\Sigma$ is compact \cite {unzush} and contains the 3-dimensional group $\Omega$, and $\Sigma$ is at most 
3-dimensional because $\Pi$ is not Clifford by assumption \cite{dim4}, hence $\Sigma^1 = \Omega$. Now $\Sigma$  is contained in 
a maximal compact subgroup ${\rm PO} (4,\BR) = \BZ_2\cdot{\rm PSO}_4\BR$ of ${\rm PGL} (4,\BR)$. 
The group $\Omega$ is the diagonal of 
${\rm PSO}_4\BR \cong \SO\times \SO$  and is normalized by $\Sigma$ since 
$\Omega = \Sigma^1$. This implies that either $\Sigma$ equals $\Omega$ itself 
or $\Sigma$ is the extension of $\Omega$ by the reflection $\varrho$ in a plane perpendicular to $Z$.
However, $\varrho$ sends every right screwed regulus in $\CC$ to a left screwed regulus, 
and this regulus is not contained in $\CC$ (no matter what the orientation is). This proves that $\Sigma = \Omega$.
As an alternative argument, one could use Lemma 3.4 of \cite{unzush} in order to rule out the possibility $\varrho \in \Sigma$.
\epf

\bprop
Let two $\So$-rotational spreads $\CC$, $\CC'$ and $\Gamma$-admissible copies $\Omega$ and $\Omega'$ of $\SO$ be given. If the 
parallelisms or oriented parallelisms $\Omega(\CC)$ and $\Omega'(\CC')$ are isomorphic, then the spreads $\CC$ and $\CC'$ 
are isomorphic. 
\eprop

\bpf
Let $\xi: \Omega(\CC) \to \Omega'(\CC')$ be an isomorphism. Since $\Omega'$ is transitive on $\Omega'(\CC')$, 
we may assume that $\xi(\CC) = \CC'$.
\epf

The converse of this implication, however plausible, is far from being true. 
The first indication for this is the fact that the additional automorphisms of the spreads discussed in Theorem \ref{dim2} do
not lead to a larger automorphism group of $\Omega(\CC)$, according to Theorem \ref{autpar}. 
In fact, we know that in the generic case, 
the automorphism group of a $\Gamma$-rotational spread is 1-dimensional, see Theorem \ref{dim2}, and then we have the following

\bthm\label{noniso}
Let $\CC$ be a $\Gamma$-rotational spread such that $\Sigma = \mathop{\rm Aut}\CC \le {\rm PGL}_4\BR$ is 1-dimensional, and let 
$\Omega \cong \SO$ be $\Gamma$-admissible (in particular, $\Omega\le \mathop{\rm AGL}_3\BR$ contains $\Phi = \Gamma^1$). 
Then there are uncountably  many conjugates $\Omega^\xi = \xi^{-1} \Omega \xi \le \mathop{\rm AGL}_3\BR$ that define pairwise non-isomorphic
parallelisms or oriented parallelisms $\Omega^\xi(\CC)$. 
\ethm

\bpf
Let $\Omega$ and $\Omega'$ both be $\Gamma$-admissible for $\CC$ and let $\Pi=\Omega(\CC)$ and $\Pi' = \Omega'(\CC)$. First we note 
that if there is an isomorphism $\rho : \Pi \to \Pi'$, then there is an isomorphism $\sigma : \Pi \to \Pi'$ contained in $\Sigma$.
This follows from transitivity of $\Omega = \mathop {\rm Aut} \Pi$. Indeed, there is $\CC' \in \Pi$ such that
$\rho (\CC') = \CC \in \Pi \cap \Pi'$, and there is $\omega \in \Omega$ such that $\omega(\CC) = \CC'$. Then
$\sigma := \rho\omega $ maps $\CC$ to itself. Now $\Sigma/\Phi$ is at most 
countable and $\Phi\le \Omega$ maps $\Pi$ to itself, 
hence $\Pi$ is isomorphic to at most countably many parallelisms $\Pi'$. The same holds for every $\Omega^\xi(\CC)$ if $\xi$ normalizes $\Phi$.

On the other hand, according to Lemma \ref{admiss}, the normalizer  of $\Phi$ in $\mathop{\rm AGL}_3\BR$ contains a 2-parameter subgroup $\Xi$ such that
the conjugates  $\Omega^\xi$ for $\xi \in \Xi$ are all $\Gamma$-admissible and pairwise distinct, namely 
the group $\Xi$ consisting of all 
affine maps $(x,y,z) \mapsto (x,y,sz+t)$ with $s, t\in \BR$ and $s>0$. (If $\Gamma = \oo$ and we want ordinary parallelisms, 
we have to set $t=0$.) 
Then the parallelisms or oriented parallelisms $\Omega^\xi(\CC)$ are also pairwise distinct, 
because every $\Omega^\xi$ is the automorphism group of $\Omega^\xi(\CC)$, by Theorem \ref{autpar}. This proves the theorem. 
\phantom{xxxx}
\epf

We remark that in the last proof we have not exhausted all possibilities of modifying $\Omega$ since we kept the hyperplane at infinity fixed.


\bibliographystyle{plain}

\begin{thebibliography}{9}

\bibitem{bettenTE}
D. Betten,
Nicht-desarguessche 4-dimensionale Ebenen,
{\em Arch. Math.} 21, 100 -- 102, 1970.

\bibitem{1fix}
D. Betten,
4-dimensionale Translationsebenen mit genau einer Fixrichtung,
{\em Geom. Dedic.} 3, 405 -- 440, 1975.

\bibitem{2fix}
D. Betten,
4-dimensionale Translationsebenen mit 7-dimensionaler Kollineationsgruppe,
{\em J. reine angew. Math.} 285, 126 -- 148, 1976.

\bibitem{komstandgp}
D.Betten,
4-dimensionale Translationsebenen mit kommutativer Standgruppe,
{\em Math. Z.} 154, 125 -- 141, 1977.

\bibitem{nonreg}
D. Betten and R. Riesinger,
Parallelisms of $\mathop{\rm PG}(3;\BR)$ composed of non-regular spreads,
{\em Aequationes Math.} 81, 227 -- 250, 2011.

\bibitem{cliff}
D. Betten and R. Riesinger,
Clifford parallelism: old and new definitions, and their use,
{\em J. Geom.} 103,  31–73, 2012.

\bibitem{coll}
D. Betten and R. Riesinger,
Collineation groups of topological parallelisms,
{\em Adv. in Geometry} 14, 175 -- 189, 2014.


\bibitem{dim3}
D. Betten and R. Riesinger,
Regular 3-dimensional parallelisms of PG(3,R), 
{\em Bull. Belg. Math. Soc. Simon Stevin} 22, 1-23, 2015.

\bibitem{zush}
D. Betten and R. L\"owen,
Compactness of the automorphism group of a topological parallelism on real projective 3-space,
{\em Results in Mathematics} 72, 1021 -- 1030, 2017.

\bibitem{kling}
W. Klingenberg,
{\em Lineare Algebra und Geometrie},
1st ed., Berlin etc.: Springer, 1984.

\bibitem{knarr}
N. Knarr,
{\em Translation Planes},
Springer Lecture Notes in Mathematics, Springer, 1995

\bibitem{dim4}
R. L\"owen,
A characterization of Clifford parallelism by automorphisms,
{\em Innovations in Incidence Geometry}, to appear; arXiv:1702.03328.

\bibitem{unzush}
R. L\"owen,
Compactness of the automorphism group of a topological parallelism on real projective 3-space: The disconnected case,
{\em Bull. Belg. Math. Soc. -- Simon Stevin}, to appear; arXiv:1710.05558.

\bibitem{CPP} 
H. Salzmann, D. Betten, T. Grundh\"ofer, H. H\"ahl, R. L\"owen, M. Stroppel, 
{\em Compact projective planes}, 
Berlin etc.: de Gruyter 1995.



\end{thebibliography}

\bigskip
\bigskip
\noindent Rainer L\"owen, Institut f\"ur Analysis und Algebra,
Technische Universit\"at Braunschweig,
Universit\"atsplatz 2,
D 38106 Braunschweig,
Germany

\end{document}